\documentclass[11pt]{article}
\usepackage[top=25mm,bottom=25mm,left=25mm,right=25mm]{geometry}

\usepackage{amsmath,amssymb,amsthm,latexsym,amscd,stmaryrd,mathrsfs,longtable,color,url}
\usepackage[dvipdfmx]{hyperref}
\usepackage{xcolor}
\hypersetup{
	colorlinks=true,
	citecolor=teal,
	linkcolor=blue,
	urlcolor=gray,
}

\newtheoremstyle{standard}%
{9pt}%
{9pt}%
{\it}
{}%
{\bfseries}%
{}
{ }%
{#3}%

\newcommand{\db}[1]{(\!({#1})\!)}

\numberwithin{equation}{section}

\newcommand{\Z}{{\mathbb Z}}
\newcommand{\Q}{{\mathbb Q}}

\newcommand{\C}{{\mathbb C}}

\newcommand{\wn}{n}

\newcommand{\wx}{x}
\newcommand{\wy}{y}

\newcommand{\ws}{s}
\renewcommand{\wr}{r}
\newcommand{\wl}{l}

\newcommand{\mN}{N}
\newcommand{\mW}{W}

\newcommand{\inda}{m}
\newcommand{\indb}{n}

\newcommand{\zhuOzero}[4]{O({#1},{#2},{#3};{#4})}

\newcommand{\unitmu}[1]{E^{(#1)}}

\newcommand{\glb}[2]{Y(#1 | #2)}
\newcommand{\ty}[4]{Y^{(#1)}_{#2}(#3 | #4)}
\newcommand{\cty}[4]{Y^{(#1)}_{#2}(#3 ; #4)}

\newcommand{\module}{M}

\newcommand{\sA}{{\mathscr A}}

\newcommand{\sM}{{\mathscr M}}

\newcommand{\sO}{{\mathscr O}}
\newcommand{\sS}{{\mathscr S}}
\newcommand{\sT}{{\mathscr T}}
\newcommand{\sU}{{\mathscr U}}
\newcommand{\sX}{{\mathscr X}}

\newcommand{\wq}{q}

\newcommand{\lu}{u}

\newcommand{\ten}{B}

\newcommand{\wz}{z}

\DeclareMathOperator{\Hom}{Hom} 
 
\DeclareMathOperator{\End}{End}

\DeclareMathOperator{\Res}{Res}

\DeclareMathOperator{\id}{id} 
\DeclareMathOperator{\Ind}{Ind}

\DeclareMathOperator{\Aut}{Aut}

\DeclareMathOperator{\Span}{Span}

\newtheorem{lemma}{Lemma}[section]
\newtheorem{theorem}[lemma]{Theorem}

\newtheorem{corollary}[lemma]{Corollary}
\theoremstyle{definition}
\newtheorem{definition}[lemma]{Definition}

\newtheorem{remark}[lemma]{Remark}

\theoremstyle{standard}

\title{A Schur-Weyl type duality for twisted weak modules over a vertex algebra}
\author{Kenichiro Tanabe\footnote{Research was partially supported by the Grant-in-aid
(No. 21K03172) for Scientific Research, JSPS.}\\\\
Faculty of Liberal Arts and Sciences\\
Tokyo City University\\
1-28-1 Tamazutsumi, Setagaya-ku,
Tokyo 158-8557\\
Japan\\
ktanabe@tcu.ac.jp}
\date{}
\begin{document}
\maketitle

\begin{abstract}
Let $V$ be a vertex algebra of countable dimension, $G$ a subgroup of $\Aut V$ 
of finite order,
$V^{G}$ the fixed point subalgebra of  $V$ under the action of  $G$,
and $\sS$ a finite $G$-stable set of inequivalent irreducible twisted 
weak $V$-modules associated with possibly different automorphisms in $G$.
We show a Schur--Weyl type duality for the actions
of $\sA_{\alpha}(G,\sS)$ and $V^G$ on the direct sum of twisted weak $V$-modules in $\sS$
where $\sA_{\alpha}(G,\sS)$ is a finite dimensional semisimple associative algebra associated with
$G,\sS$, and a $2$-cocycle $\alpha$ naturally determined by the $G$-action on $\sS$. 
It follows as a natural consequence of the result that for any $g\in G$ every irreducible
$g$-twisted weak $V$-module is a completely reducible weak $V^G$-module.
\end{abstract}



\section{\label{section:introduction}Introduction}
Let $V$ be a vertex algebra (cf. \cite{B1986, LL}) and $G$ a subgroup of $\Aut V$ of finite order.
We denote by $V^{G}$ the fixed point subalgebra of  $V$ under the action of  $G$: $V^{G}=\{a\in V\ |\ g(a)=a\mbox{ for all }g\in G\}$.
The fixed point subalgebras play an important role in the study of vertex algebras, particularly in the 
construction of interesting examples with the moonshine vertex algebra $V^{\natural}$ as a representative example \cite{B1986,FLM}.
One of the main problems with $V^G$ is  describing the $V^G$-modules in terms of $V$ and $G$.
In line with the various ideas proposed in \cite{DVVV1989},
for a simple vertex operator algebra $V$,
Dong and Mason initiated a systematic study of representations of $V^G$
and showed a Schur-Weyl type duality  for $(G,V^{G})$ when $G$ is a solvable group \cite{DM1997}.
Here, it is worth mentioning that as one consequence of the Schur-Weyl type duality, 
 a Galois correspondence in simple vertex operator algebras was established in \cite{DM1997} when $G$ is either abelian or dihedral,
and in \cite{HMT1999} in the general case.
The Schur-Weyl type duality was extended in \cite{DM1996, DY2002, MT2004, Yamskphd}
using the Zhu algebra in \cite{Z1996} and its generalizations in \cite{DLM1998t,DLM1998tn,DLM1998v, MT2004}.  
 In connection with this problem, twisted $V$-modules have been studied systematically (cf. \cite{DLM1998t,FLM,Lepowsky1985,Li1996})
for the following reasons: for $g\in G$, every $g$-twisted $V$-module becomes a $V^G$-module
and, moreover, it is conjectured that under some conditions on $V$,
every irreducible $V^G$-module is contained in some irreducible $g$-twisted 
$V$-module for some $g\in G$ (cf. \cite{DVVV1989}).
Let $\sS$ be a finite $G$-stable set (see the explanation just after \eqref{eqn:actionYmodulecdotsigma} for the definition) of inequivalent 
irreducible 
twisted $V$-modules
associated with possibly different automorphisms in $G$.
A finite dimensional semisimple associative $\C$-algebra $\sA_{\alpha}(G,\sS)$ (see \eqref{eq:algebraAalphaGS} for the definition) associated to $G$, $\sS$, and a
 $2$-cocycle $\alpha$ naturally determined by the $G$-action on $\sS$ was constructed in \cite{DY2002}.
The algebra $\sA_{\alpha}(G,\sS)$ acts on the direct sum $\sM=\oplus_{M\in\sS}M$
and the actions of $\sA_{\alpha}(G,\sS)$ and  $V^G$ commute with each other. 
The Schur-Weyl type duality  above  is about a decomposition of $\sM$ under the action of $\sA_{\alpha}(G,\sS)\otimes_{\C}V^G$. 
In \cite{DM1997, DY2002, Yamskphd} they studied the case that $\sS$ consists of $g$-twisted $V$-modules associated with a single automorphism $g\in G$,
where we note that in this case, $\sM$ is also a $g$-twisted $V$-module.
If this is not the case,
then $\sM$ is not a twisted $V$-module in general since
the definition of a twisted $V$-module depends on each automorphism of $V$
and hence a direct sum of twisted $V$-modules associated with different automorphisms is not always a twisted $V$-module.
This fact makes it difficult to investigate $\sM$.
To overcome the difficulty, Miyamoto and the author introduced the $G$-twisted 
Zhu algebra which allows us to 
study twisted $V$-modules in  a  unified way and we showed a Schur-Weyl type duality 
for any finite $G$-stable set $\sS$ of inequivalent irreducible twisted $V$-modules in \cite[Theorem 2]{MT2004}.

The purpose of this paper is to generalize the Schur-Weyl type duality in \cite[Theorem 2]{MT2004} to twisted weak modules for a vertex algebra. 
The settings are all the same as above if we replace the terms  \lq\lq vertex operator algebra\rq\rq\ and \lq\lq(twisted) $V$-module\rq\rq\ by \lq\lq vertex algebra\rq\rq\ and \lq\lq (twisted) weak $V$-module\rq\rq,\ respectively.
Here we mention that the weak modules for $V_{L}^{+}$, a fixed point subalgebra of the lattice vertex algebra $V_{L}$ where 
$L$ is a non-degenerate even lattice, are classified in \cite{Tanabe2021-1} and every weak $V_{L}^{+}$-module is a submodule of some twisted weak $V_{L}^{+}$-module.
The same result is shown for the weak modules with some conditions for $M(1)^{+}$,
a fixed point subalgebra of the Heisenberg vertex operator algebra $M(1)$, in \cite{Tanabe2017}.
These examples lead us to expect a Schur-Weyl type duality holds for twisted 
weak $V$-modules.
For an automorphism $g \in \Aut V$ of finite order, we denote 
by $\sT_{g}$ the set of all irreducible $g$-twisted weak $V$-modules
and we define  $\sT_{G}=\cup_{g\in G}\sT_{g}$.
Let $\sS$ be a finite $G$-stable subset of $\sT_{G}$ consisting of inequivalent 
irreducible 
twisted weak $V$-modules
and let $\sM=\oplus_{M\in\sS}M$.
Under the assumption that $V$ is a simple vertex algebra of countable dimension, a Schur-Weyl type duality for $(\sA_{\alpha}(G,\sS),V^G)$ was studied in \cite{ALPY2019, ALPY2022, DRY2023}.
There it is assumed that $\sS$ consists of $g$-twisted weak $V$-modules associated with a single automorphism $g\in G$. 
If this is not the case, as in the case of twisted $V$-modules, $\sM$ is not a twisted weak $V$-module in general.
Moreover, we can not apply the $G$-twisted Zhu algebra to twisted weak $V$-modules
since we do not assume any grading in the definition of a twisted weak $V$-module. 
To overcome the difficulty,
instead of the $G$-twisted Zhu algebra, we shall use weak $(V,T)$-modules (see Definition \ref{definition:def-vt} for the definition) 
defined for any positive integer $T$, introduced by the author \cite{Tanabe2015}
as a generalization of the notion of a twisted weak module so that it is closed under taking \textcolor{black}{direct sums}.
For any $g\in \Aut V$ of finite order and any $T$ that is a multiple of the order of $g$, every $g$-twisted weak $V$-module becomes a weak $(V,T)$-module.
The notion of a weak $(V,T)$-module  allows us to study twisted weak $V$-modules in a unified way
(Lemma \ref{lemma:twisted-(V,T)}).
We state the main result of this paper, which is a  Schur-Weyl type duality for $(\sA_{\alpha}(G,\sS),V^G)$
in the general case:
\begin{theorem}
\label{theorem:main}
Let $V$ be a simple vertex algebra of countable dimension,
$G$ a subgroup of $\Aut V$ of finite order,
$\sS$ a finite $G$-stable set of inequivalent weak $(V,T)$-modules,
and $\sM=\oplus_{\module\in\sS}\module$.
Under the action of $\sA_{\alpha}(G,\sS)\otimes_{\C}\textcolor{black}{V^G}$, the space $\sM$ decomposes into a direct sum
\begin{align}
\sM&=\bigoplus\limits_{j\in J, \lambda\in\Lambda_{j}}W_{\lambda}^{j}\otimes_{\C}M_{\lambda}^{j}.
\end{align}
Moreover
\begin{enumerate}
\item For any $j\in J$ and $\lambda\in\Lambda_j$, $M^{j}_{\lambda}$ is non-zero and is an irreducible weak $(V^{G},T)$-module.
\item For $j_1,j_2\in J$ and $\lambda_1\in\Lambda_{j_1},\lambda_2\in\Lambda_{j_2}$,
$M^{j_1}_{\lambda_1}$ and $M^{j_2}_{\lambda_2}$ are isomorphic as weak $(V^G,T)$-modules 
if and only if $(j_1,\lambda_1)=(j_2,\lambda_2)$.   
\end{enumerate}
\end{theorem}
\noindent{}In the theorem, $\{W^{j}_{\lambda}\ |\ j\in J, \lambda\in\Lambda_{j}\}$
is a complete set of representatives of isomorphism
classes of irreducible $\sA_{\alpha}(G,\sS)$-modules
and $M^{j}_{\lambda}$ is the multiplicity space of $W^{j}_{\lambda}$ in $\sM$ for $j\in J$ and $\lambda\in \Lambda_{j}$,
on  which $V^G$ acts naturally (see  \eqref{eq:Mlambda=HomC}, \eqref{eq:definitionWjlambda},
and the explanation just before \eqref{eq:GjGmoduelj} for the precise definitions of
$W^{j}_{\lambda}$ and  $M^{j}_{\lambda}$).

We note that if $\sS$ is a $G$-stable subset of
$\sT_{G}$ consisting of inequivalent irreducible twisted weak $V$-modules,
then $\sM=\oplus_{M\in\sS}M$ becomes a weak $V^G$-module.
Hence, in this case, from the definition of  a weak $(V,T)$-module (Definition \ref{definition:def-vt}) we can replace the term \lq\lq weak $(V^G,T)$-module\rq\rq\ by  \lq\lq weak $V^G$-module\rq\rq\ in Theorem \ref{theorem:main}.
Setting $T=|G|$, we have the following result as a special case of Theorem \ref{theorem:main}:
\begin{corollary}
\label{corollary:stabletwisted}
Let $V$ be a simple vertex algebra of countable dimension,
$G$ a subgroup of $\Aut V$ of finite order,
$\sS$ a finite $G$-stable subset of $\sT_{G}$ consisting of inequivalent 
irreducible twisted weak $V$-modules,
and $\sM=\oplus_{\module\in\sS}\module$.
Under the action of $\sA_{\alpha}(G,\sS)\otimes_{\C}\textcolor{black}{V^G}$, the space $\sM$ decomposes into a direct sum
\begin{align}
\sM&=\bigoplus\limits_{j\in J, \lambda\in\Lambda_{j}}W_{\lambda}^{j}\otimes_{\C}M_{\lambda}^{j}.
\end{align}
Moreover
\begin{enumerate}
\item For any $j\in J$ and $\lambda\in\Lambda_{j}$, $M^{j}_{\lambda}$ is non-zero and is an irreducible weak $V^{G}$-module.
\item For $j_1,j_2\in J$ and $\lambda_1\in\Lambda_{j_1},\lambda_2\in\Lambda_{j_2}$,
$M^{j_1}_{\lambda_1}$ and $M^{j_2}_{\lambda_2}$ are isomorphic as weak $V^G$-modules 
if and only if $(j_1,\lambda_1)=(j_2,\lambda_2)$.   
\end{enumerate}
In particular, for any $g\in G$ every irreducible $g$-twisted weak $V$-module is a completely reducible weak $V^G$-module.
\end{corollary}

We briefly introduce the basic idea.
The key result is Lemma \ref{lemma:associative-action} below, which is a kind of \lq\lq associativity\rq\rq\ of 
the action of $V$ on a weak $(V,T)$-module $\module$ in the following sense:
for $a,b\in V$, $u\in\module$, and  $\inda,\indb\in(1/T)\Z$,
\begin{align}
a_{\inda}b_{\indb}u&\in \Span_{\C}\{(a_ib)_{\inda+\indb-i}\lu\in\module\ |\ i\in\Z\}.
\end{align}
The result is well-known when $\module$ is a (twisted) weak $V$-module (cf. \cite[Proposition 4.5.7]{LL}).
However, in the general case, the result is far from trivial 
due to the complexity of the structure of a weak $(V,T)$-module.
In the proof, we use the methods developed in \cite{MT2004} and \cite{Tanabe2015}.
Once Lemma \ref{lemma:associative-action} is established, the argument in \cite[Section 7]{DRY2023} works well
even for 
weak $(V,T)$-modules,
which leads Theorem \ref{theorem:main}.

The organization of the paper is as follows.
\textcolor{black}{In Section \ref{section:preliminary}} we recall some basic properties of 
weak $(V,T)$-modules and show some results.
In Section \ref{section:proof} after recalling 
some properties of the algebra $\sA_{\alpha}(G,\sS)$, we show the main theorem.

\section{\label{section:preliminary}Preliminary}
Throughout this paper, $\Z$ denotes the set of all integers,
$T$ is a fixed positive integer, and
$(V,Y,{\mathbf 1})$ is a vertex algebra.
Recall that $V$ is the underlying vector space, 
$Y(\mbox{ },\wx)$ is a linear map from $V\otimes_{\C}V$ to $V\db{x}$,
and ${\mathbf 1}$ is the vacuum vector.
For $g\in\Aut V$, $|g|$ denotes the order of $g$.
For $g\in\Aut V$ or  finite order, a positive integer $n$ that is a multiple of $|g|$,
and $\wr=0,1,\ldots,\wn-1$,
we define $V^{(g,\wr)}=\{a\in V\ |\ ga=e^{-2\pi\sqrt{-1}\wr/n}a\}$.
For $a\in V$, $a^{(g,\wr)}$ denotes the $\wr$th component of $a$ in the decomposition $V=\oplus_{\wr=0}^{n-1}V^{(g,\wr)}$,
that is, $a=\sum_{\wr=0}^{\wn-1}a^{(g,\wr)}$, $a^{(g,\wr)}\in V^{(g,\wr)}$.
For $i,j\in\Z$, we define
\begin{align}
\Z_{\leq i}&=\{k\in\Z\ |\ k\leq i\},\nonumber\\
\Z_{\geq i}&=\{k\in\Z\ |\ k\geq i\},\nonumber\\
\C[\wz,\wz^{-1}]_{\leq i}&=\Span_{\C}\{\wz^k\ |\ k\leq i\},\mbox{ 
and}\nonumber\\
\C[\wz,\wz^{-1}]_{i,j}&=\Span_{\C}\{\wz^k\ |\ i\leq k\leq j\}.
\end{align}

\subsection{\label{section:subspace}
Subspaces of $\C[\wz,\wz^{-1}]$
}
For $\gamma,\wq\in\Z$ and $\alpha\in\Q$, we denote 
by $O(\gamma,\alpha,\wq;\wz)$ the subspace of $\C[\wz,\wz^{-1}]$
spanned by
\begin{align}
\label{eq:res-x-q}
\Res_{\wx}\big((1+\wx)^{\alpha+1}\wx^{\wq-j}\sum_{i\in\Z_{\leq \gamma}}\wz^i\wx^{-i-1}\big)
&=
\sum\limits_{i=0}^{\gamma-\wq+j}\binom{\alpha+1}{i}\wz^{i+\wq-j},\quad j=0,1,\ldots
\end{align}
and $\wz^i, i\in\Z_{\geq \gamma+1}$, which  will be used in Lemmas \ref{lemma:vanish-Oz} and \ref{lemma:associative-action}.
The subspace $\zhuOzero{\gamma}{\alpha}{\wq}{\wz}$ in this paper is written as $\zhuOzero{\gamma}{\alpha+1}{\wq}{\wz}$ in \cite{Tanabe2015}.
We note that 
\begin{align}
\label{eq:osubset}
O(\gamma,\alpha,q;\wz)\subset O(\gamma,\alpha,q+1;\wz)\subset\cdots.
\end{align}

The following result holds as a special case of \cite[Lemma 3.2]{Tanabe2015} 
and its proof. 
\begin{lemma}\label{lemma:iso-oplus}
For $\gamma,\wq\in\Z$ and $\alpha\in \Q$, 
the diagonal map $\C[\wz,\wz^{-1}]\ni f\mapsto (f,\ldots,f)\in \C[\wz,\wz^{-1}]^{\oplus T}$
induces an isomorphism 
\begin{align}
\label{eq:C[wzwz-1]/bigcap}
\C[\wz,\wz^{-1}]/\bigcap_{\ws=0}^{T-1}O(\gamma,\alpha+\frac{\ws}{T},\wq;\wz)
&\overset{\cong}{\rightarrow} \bigoplus_{{\ws}=0}^{T-1}\C[\wz,\wz^{-1}]/O(\gamma,\alpha+\frac{\ws}{T},\wq;\wz)
\end{align}
as vector spaces.
Moreover, each representative element of the quotient space on the left-hand side of \eqref{eq:C[wzwz-1]/bigcap}
can be taken uniquely from $\C[z,z^{-1}]_{\gamma+1-T(\gamma-\wq),\gamma}$.
\end{lemma}

For $\wr=0,\dots,T-1$ and $i\in\Z_{\leq \gamma}$, by
Lemma \ref{lemma:iso-oplus}
there exists a unique $\unitmu{\wr}(\gamma,\alpha,\wq,i;\wz)\in 
\C[z,z^{-1}]_{\gamma+1-T(\gamma-\wq),\gamma}$
such that
\begin{align}
\label{eq:unit-s}
\unitmu{\wr}(\gamma, \alpha,\wq,i;\wz)
&\equiv \delta_{r,s}z^i\pmod{
\zhuOzero{\gamma}{\alpha+\frac{\ws}{T}}{\wq}{\wz}},\qquad s=0,1,\ldots,T-1.
\end{align}
We also define 
$\unitmu{\wr}(\gamma,\alpha,\wq,i;\wz)=0\mbox{ for $i\in\Z_{\geq \gamma+1}$}$ for convenience.

\subsection{$(V,T)$-modules
\label{section:vt}
}
In this subsection, we recall the definition of  a weak $(V,T)$-module
and its properties from \cite{Tanabe2015}.
We also show some results about weak $(V,T)$-modules.
For a vector space $M$ over $\C$,
we define three linear injective maps
\begin{align*}
\iota_{\wx,y} : & M[[\wx^{1/T},y^{1/T}]][\wx^{-1/T},y^{-1/T},(\wx-y)^{-1}]\rightarrow M\db{\wx^{1/T}}\db{y^{1/T}},\\
\iota_{y,\wx} : &M[[\wx^{1/T},y^{1/T}]][\wx^{-1/T},y^{-1/T},(\wx-y)^{-1}]\rightarrow M\db{y^{1/T}}\db{\wx^{1/T}},\\
\iota_{\wx,\wx-y} : & M[[\wx^{1/T},y^{1/T}]][\wx^{-1/T},y^{-1/T},(\wx-y)^{-1}]\rightarrow M\db{y^{1/T}}\db{\wx-y}
\end{align*}
by 
\begin{align*}
\iota_{\wx,y}
f&=\sum_{j,k,l}u_{j,k,l}\sum_{i=0}^{\infty}\binom{l}{i}(-1)^{i}\wx^{j+l-i}y^{k+i},\\
\iota_{y,\wx}
f&=\sum_{j,k,l}u_{j,k,l}\sum_{i=0}^{\infty}\binom{l}{i}(-1)^{l-i}y^{k+l-i}\wx^{j+i},\\
\iota_{y,\wx-y}
f&=\sum_{j,k,l}u_{j,k,l}\sum_{i=0}^{\infty}\binom{j}{i}y^{k+j-i}(\wx-y)^{l+i}
\end{align*}
for $f=\sum_{j,k,l}u_{j,k,l}
\wx^jy^k(\wx-y)^l\in M[[\wx^{1/T},y^{1/T}]][\wx^{-1/T},y^{-1/T},(\wx-y)^{-1}],
u_{j,k,l}\in M$.

Now we write down the definition of a {\it weak $(V,T)$-module}, which  is called a $(V,T)$-module in \cite[Definition 2.1]{Tanabe2015}.
\begin{definition}\label{definition:def-vt}
A {\it weak $(V,T)$-module} $\module$ is a vector space over $\C$ equipped with a linear map
\begin{align}
\label{eq:inter-form}
Y_{\module}(\ , x) : V\otimes_{\C}\module&\rightarrow \module\db{x}\nonumber\\
a\otimes u&\mapsto  Y_{\module}(a, x)\lu=\sum_{n\in(1/T)\Z}a_{n}\lu x^{-n-1}
\end{align}
such that the following conditions are satisfied:
\begin{enumerate}
\item $Y_{M}({\bf 1},\wx)=\id_{M}$.
\item For $a,b\in V$ and $\lu\in M$,
there exists $\glb{a,b,\lu}{\wx,\wy}\in M[[\wx^{1/T},y^{1/T}]][\wx^{-1/T},y^{-1/T},(\wx-y)^{-1}]$
such that 
\begin{align*}
\iota_{\wx,y}\glb{a,b,\lu}{\wx,\wy}&=Y_{M}(a,\wx)Y_{M}(b,y)\lu,\\
\iota_{y,\wx}\glb{a,b,\lu}{\wx,\wy}&=Y_{M}(b,y)Y_{M}(a,\wx)\lu,\quad\mbox{and }\\
\iota_{y,\wx-y}\glb{a,b,\lu}{\wx,\wy}&=Y_{M}(Y(a,\wx-y)b,y)\lu.
\end{align*}
\end{enumerate}
\end{definition}

For a weak $(V,T)$-module $M$, 
a subspace $N$ of $M$ is called {\it weak $(V,T)$-submodule} of $M$
if $(N,Y_M|_{N})$ is a weak $(V,T)$-module,
where $Y_M|_{N}$ is the restriction of $Y_{M}$ to $N$.
A non-zero weak $(V,T)$-module $M$ is called {\it irreducible} or {\it simple} 
if there is no submodule of $M$ except $0$ and $M$ itself.
For weak $(V,T)$-modules $\module$ and $\mN$, a linear map $f : \module\rightarrow \mN$
is called a {\it homomorphism} if $f(a_iu)=a_if(u)$ for all $a\in V$, $u\in \module$, and $i\in (1/T)\Z$.
For weak $(V,T)$-modules $\module$ and $\mN$, we call $\module$ is \textcolor{black}{\it isomorphic} to $\mN$ as weak $(V,T)$-modules,
which we write as $\module\cong\mN$,
if there exist homomorphisms $f : \module\rightarrow \mN$ and $g : \mN\rightarrow\module$ such that
$f\circ g=\id_{\mN}$ and 
$g\circ f=\id_{\module}$.
For a submodule $N$ of a weak $(V,T)$-module $M$,
the quotient space $M/N$ is a weak $(V,T)$-module.
For a set of weak $(V,T)$-modules $\{M_{i}\}_{i\in I}$,
the direct sum $\oplus_{i\in I}M_i$ is  a weak $(V,T)$-module.

Let $M$ be a weak $(V,T)$-module.
For $a\in V$ and $s=0,1,\ldots,T-1$, we define 
\begin{align}
\label{eq:Ys}
Y^{s}_{M}(a,\wx)&=\sum_{\begin{subarray}{c}i\in s/T+\Z\end{subarray}}
a_{i}\wx^{-i-1}.
\end{align}
Clearly, $\sum_{s=0}^{T-1}Y^{s}_{M}(a,\wx)=Y_{\module}(a,\wx)$ holds.
For $a,b\in V$, $\lu\in\module$, and $s=0,1,\ldots,T-1$,  writing
\begin{align}
Y_{\module}(a,b,\lu | x_1, x_2)&=\sum_{i,j\in(1/T)\Z}\sum_{k\in\Z}v_{ijk}x_1^{i}x_2^{j}(x_1-x_2)^{k},\ v_{ijk}\in \module,
\end{align}
we define
\begin{align}
Y^{(s)}_{\module}(a,b | x_2, x_1-x_2)(u)&=\iota_{x_2,x_1-x_2}\big(\sum_{i\in -s/T+\Z}\sum_{j\in (1/T)\Z}\sum_{k\in\Z}v_{ijk}x_1^{i}x_2^{j}(x_1-x_2)^{k}\big).
\end{align}
By the definition,
\begin{align}
\label{eq:sum-Y}
\sum_{s=0}^{T-1}\ty{s}{M}{a,b}{\wx_2,\wx_0}(\lu)&=Y_{M}(Y(a,\wx_0)b,\wx_2)\lu.
\end{align}
For $q\in\Z$, $M\db{\wx_2^{1/T}}\db{\wx_0}_{\geq q}$ denotes
the set of all elements in $M\db{\wx_2^{1/T}}\db{\wx_0}$ of the form
$\sum_{\begin{subarray}{c}i\in(1/T)\Z\\j\in\Z_{\geq q}\end{subarray}}
u_{ij}\wx_2^{i}\wx_0^{j},\ u_{ij}\in\module$.
By \cite[Remark 2.6]{Tanabe2015}, if $Y_{M}(Y(a,\wx_0)b,\wx_2)\lu\in M\db{\wx_2^{1/T}}\db{\wx_0}_{\geq q}$,
then so is $Y^{(s)}_{\module}(a,b | x_2, x_0)(u)$ for any $s=0,1,\ldots,T-1$.
We write the expansion of $\ty{s}{M}{a,b}{\wx_2,\wx_0}$ as
\begin{align}
\label{eq:tysMabwx2wx0}
\ty{s}{M}{a,b}{\wx_2,\wx_0}&=\sum_{i\in(1/T)\Z}\sum_{j\in\Z}\cty{s}{M}{a,b}{i,j}\wx_2^{-i-1}\wx_0^{-j-1},\ \cty{s}{M}{a,b}{i,j}\in \End_{\C}M.
\end{align}
We recall the following result from \cite[(2.16)]{Tanabe2015},
which is essentially the same as the Borcherds identity:
for $a,b\in V$, $\lu\in\module$, and $j,k\in(1/T)\Z, l\in\Z$,
\begin{align}
\label{eq:Borcherds-Coeff}
\sum_{i=0}^{\infty}\binom{j}{i}\cty{s}{M}{a,b}{j+k-i,l+i}(\lu)
&=
\sum_{i=0}^{\infty}\binom{l}{i}(-1)^{i}(a_{l+j-i}b_{k+i}+(-1)^{l+1}b_{l+k-i}a_{j+i})\lu.
\end{align}

\textcolor{black}{We recall} the relationship between the notion of a weak $(V,T)$-module and that of a (twisted) weak $V$-module
 discussed in \cite{Tanabe2015}.
Before that, we caution that  if $V$ is a vertex operator algebra, then 
the notion of a module for $V$ viewed as a vertex algebra is different from the notion of a module for $V$ viewed as a vertex operator algebra (cf. \cite[Definitions 4.1.1 and 4.1.6]{LL}).
To avoid confusion, throughout this paper, we refer to a module for a vertex algebra defined in \cite[Definition 4.1.1]{LL} 
as a {\it weak module} (cf. \cite[Definition 2.1]{Tanabe2021-1}).
The reason why we use the terminology \lq\lq weak module\rq\rq\ is that when $V$ is a vertex operator algebra, a module for $V$ viewed as a vertex algebra is called 
a weak $V$-module (cf. \cite[Definition 2.3]{ABD2004}, \cite[p.150]{DLM1997}, and \cite[p.157]{Li1996}).
We also apply this rule for twisted modules and twisted weak modules.
By \cite[Lemma 2.4]{Tanabe2015}, we find that the following definitions of 
a (twisted) weak $V$-module are the same as the standard ones (\cite[Definition 2.1]{Tanabe2021-1} and \cite[Definition 2.6]{Li1996t}):
\begin{definition}\label{definition:def-weak-twisted-weak}
\begin{enumerate}
\item A weak $(V,1)$-module is called a {\it weak $V$-module}.
\item
Let $g\in \Aut V$ of finite order.
A weak $(V,|g|)$-module $\module$ is called a {\it $g$-twisted weak $V$-module} if 
$Y_{\module}(a,x)=Y^{s}_{\module}(a,x)$ for all $s=0,1,\ldots,|g|-1$ and $a\in V^{(g,s)}$.  
\end{enumerate}
\end{definition}

The following result implies that the notion of a weak $(V,T)$-module  allows us to study twisted weak $V$-modules in a unified way.
\begin{lemma}
\label{lemma:twisted-(V,T)}
Let $g,h\in\Aut V$ of finite order and $T$ a common multiple of $|g|$ and $|h|$.
Let $\module$ be a $g$-twisted weak $V$-module and 
$\mN$ an $h$-twisted weak $V$-module. Then, $\module \cong \mN$ as weak $(V,T)$-modules
if and only if $g=h$ and $\module\cong\mN$ as $g$-twisted weak modules.
\end{lemma}
\begin{proof}
Assume $\module\cong\mN$ as weak $(V,T)$-modules.
\textcolor{black}{We} note that $V^{(g,rT/|g|)}\neq 0$ for all $r=0,1,\ldots,|g|-1$ since $V$ is simple.
For any $r=0,1,\ldots,|g|-1$ and $0\neq a\in V^{(g,rT/|g|)}$,
since $V$ is simple, $0\neq Y_{\module}(a,x)=Y_{\module}^{rT/|g|}(a,x)$ and hence $Y_{\mN}(a,x)=Y^{rT/|g|}_{\mN}(a,x)$.
Thus, $g=h$.
The rest of the  proof is clear.
\end{proof}

For a vertex algebra $V$, we take the tensor algebra $\ten(V)=T(V\otimes_{\C} \C[\zeta^{1/T},\zeta^{-1/T}])$ where $\zeta^{1/T}$ is a variable.
For a weak $(V,T)$-module $\module$, we regard $\module$ as a $\ten(V)$-module in the following way: for 
$a\in V,  i\in(1/T)\Z$, and $\lu\in \module$, we define
\begin{align}
(a\otimes \zeta^{i})\cdot \lu&=a_{i}u.
\end{align}
We note that if $V$ is  of countable dimension, then so is $B(V)$. 
The following result clearly holds:
\begin{lemma}
\label{lemma:tensor}
Let $V$ be a vertex algebra and let $\module$, $\mN$ be weak $(V,T)$-modules.
Then
\begin{enumerate}
\item $\module$ is an irreducible  weak $(V,T)$-module if and only if $\module$ is an irreducible  weak $\ten(V)$-module.
\item $\module \cong \mN$ as  weak $(V,T)$-modules  if and only if $\module \cong \mN$ as $\ten(V)$-modules.
\item (Schur's lemma) If $V$ is of countable dimension and $\module$ is irreducible, then $\End_{V}\module=\End_{\ten(V)}\module\cong \C$.
\end{enumerate}
\end{lemma}
Let $\module$ be a weak $(V,T)$-module.
For $a\in V$ and $\lu\in \module$, we define  $\epsilon(a,\lu)\in(1/T)\Z\cup\{-\infty\}$ by
\begin{align}
\label{eqn:max-vanish0}
a_{\epsilon(a,\lu)}\lu&\neq 0\mbox{ and }a_{i}\lu
=0\mbox{ for all }i>\epsilon(a,\lu)
\end{align}
if $Y_{\module}(a,x)\lu\neq 0$ and $\epsilon(a,\lu)=-\infty$
if $Y_{\module}(a,x)\lu= 0$.

For the next results, we prepare the following symbols.
Let $M$ be a weak $(V,T)$-module,
$a,b\in V$, $u\in \module$, $\inda,\indb\in(1/T)\Z$, and $\alpha,\beta,\gamma \in\Z$ such that 
$\alpha\geq \epsilon(a,u)$, $\beta\geq \epsilon(b,u)$, and $\gamma\geq\epsilon(a,b)$.
We write $\inda=\wl_1+\wr_1/T$ and $\indb=\wl_2+\wr_2/T$ with $\wl_1,\wl_2\in\Z$ and $0\leq \wr_1,\wr_2<T$.
For $s=0,1,\ldots,T-1$,
we define a linear map $Z^{(\ws)}_{\inda,\indb}(a,b,\lu;-) : 
\C[\wz,\wz^{-1}]\rightarrow M$ by
\begin{align}
\label{eq:zmnm}
Z^{(\ws)}_{\inda,\indb}(a,b,\lu;\wz^i)&=\cty{s}{M}{a,b}{\inda+\indb-i,i}(\lu),\ i\in\Z.
\end{align}
Here $\cty{s}{M}{a,b}{\inda+\indb-i,i}(\lu)$ is defined in 
\eqref{eq:tysMabwx2wx0}.
\begin{lemma}
\label{lemma:vanish-Oz}
For $s=0,1,\ldots,T-1$ and $p(z)\in\zhuOzero{\gamma}{\alpha+s/T}{\wl_1-\alpha+\wl_2-\beta-2}{\wz}$, $Z^{(\ws)}_{\inda,\indb}(a,b,\lu;p(\wz))=0$.
\end{lemma}
\begin{proof}
If $i\geq \gamma+1$, then $Z^{(\ws)}_{\inda,\indb}(a,b,\lu;z^i)=0$ by 
\eqref{eq:zmnm} and the explanation just after \eqref{eq:sum-Y}.
By  \eqref{eq:res-x-q}, \eqref{eq:Borcherds-Coeff}, \eqref{eqn:max-vanish0}, and \eqref{eq:zmnm},
for any $j=0,1,\ldots$,
\begin{align}
&Z^{(\ws)}_{\inda,\indb}(a,b,\lu;\sum_{i=0}^{\gamma}
\binom{\alpha+s/T+1}{i}\wz^{i+\wl_1-\alpha+\wl_2-\beta-2-j})\nonumber\\
&=\sum_{i=0}^{\gamma}\binom{\alpha+s/T+1}{i}\cty{s}{M}{a,b}{\inda+\indb-i-\wl_1+\alpha-\wl_2+\beta+2+j,i+\wl_1-\alpha+\wl_2-\beta-2-j}(\lu)
\nonumber\\
&=\sum_{i=0}^{\infty}\binom{\alpha+s/T+1}{i}\nonumber\\
&\quad{}\times\cty{s}{M}{a,b}{(\alpha+s/T+1)+(\beta+(\wr_1+\wr_2-s)/T+1+j)-i,\wl_1-\alpha+\wl_2-\beta-2-j+i}(\lu)
\nonumber\\
&=\sum_{i=0}^{\infty}\binom{\wl_1-\alpha+\wl_2-\beta-2-j}{i}(-1)^{i}\big(a_{s/T+\wl_1+\wl_2-\beta-j-1-i}b_{\beta+(\wr_1+\wr_2-s)/T+1+j+i}\nonumber\\
&\qquad{}
+(-1)^{\wl_1-\alpha+\wl_2-\beta-j-1}b_{(\wr_1+\wr_2-s)/T+\wl_1-\alpha+\wl_2-1-i}a_{\alpha+s/T+1+i}\big)\lu
\nonumber\\
&=0.
\end{align}
Here we have used $-1<(\wr_1+\wr_2-s)/T$. The proof is complete.
\end{proof}
The following is the key result of this paper.
\begin{lemma}
\label{lemma:associative-action}
Let $\module$ be a weak $(V,T)$-module.
For $a,b\in V, u\in \module$, and  $\inda,\indb\in(1/T)\Z$,
\begin{align}
a_{\inda}b_{\indb}u&\in \Span_{\C}\{(a_ib)_{\inda+\indb-i}\lu\in\module\ |\ i\in\Z\}.
\end{align}
\end{lemma}
\begin{proof}
Let $\alpha,\beta,\gamma \in\Z$ such that $\alpha\geq \epsilon(a,u)$, $\beta\geq \epsilon(b,u)$, and $\gamma\geq\epsilon(a,b)$.
We write $\inda=\wl_1+\wr_1/T$ and $\indb=\wl_2+\wr_2/T$ with $\wl_1,\wl_2\in\Z$ and $0\leq \wr_1,\wr_2<T$.
We take $\mu,\nu\in (1/T)\Z$ so that $\epsilon(a,u)\leq \mu, \epsilon(b,u)\leq \nu,$ and $\inda\equiv \mu, \indb\equiv \nu\pmod{\Z}$.
We define a Laurent polynomial $\Phi_{\inda, \indb}(\wz)\in\C[z,z^{-1}]$ by
\begin{align}
\label{eq:def-mul-x}
&\Phi_{\inda, \indb}(\wz)\nonumber\\
&=\sum_{j=0}^{\nu-\indb}
\binom{\inda-\mu-1}{j}
\Res_{x}(1+x)^{\mu+1}x^{\inda-\mu-1-j}\sum_{i\in\Z_{\leq\gamma}}
\unitmu{\wr_1}(\gamma,\alpha,\wl_1-\alpha+\wl_2-\beta,i;\wz)\wx^{-i-1}
\nonumber\\
&=\sum_{j=0}^{\nu-\indb}
\binom{\inda-\mu-1}{j}\sum_{i=0}^{\gamma-(\inda-\mu-1-j)}\binom{\mu+1}{i}
\unitmu{\wr_1}(\gamma,\alpha,\wl_1-\alpha+\wl_2-\beta,i+\inda-\mu-1-j;\wz)
\nonumber\\
&\in\C[z,z^{-1}]_{\gamma+1-T(\gamma-\wl_1-\wl_2+\alpha+\beta),\gamma}.
\end{align}
For any $i\in\Z$ and $0\leq \ws\leq T-1$ with $\ws\neq \wr_1$,
by \eqref{eq:unit-s},
\begin{align}
\unitmu{\wr_1}(\gamma,\alpha,\wl_1-\alpha+\wl_2-\beta,i;\wz)&\equiv 0\pmod{\zhuOzero{\gamma}{\alpha+\dfrac{s}{T}}{\wl_1-\alpha+\wl_2-\beta}{\wz}}
\end{align}
and hence by Lemma \ref{lemma:vanish-Oz}, $Z^{(\ws)}_{\inda,\indb}(a,b,\lu;\Phi_{\inda,\indb}(z))=0$.
By Lemma \ref{lemma:vanish-Oz}, we also have
\begin{align}
\Phi_{\inda, \indb}(\wz)
&=\sum_{j=0}^{\nu-\indb}
\binom{\inda-\mu-1}{j}\sum_{i=0}^{\gamma-(\inda-\mu-1-j)}\binom{\mu+1}{i}
z^{i+\inda-\mu-1-j}
\pmod{\zhuOzero{\gamma}{\alpha+\dfrac{r_1}{T}}{\wl_1-\alpha+\wl_2-\beta}{\wz}}
\end{align}
and hence by \eqref{eq:Borcherds-Coeff},
\begin{align}
&Z^{(\wr_1)}_{\inda,\indb}(a,b,\lu;\Phi_{\inda,\indb}(z))\nonumber\\
&=
\sum_{j=0}^{\nu-\indb}
\binom{\inda-\mu-1}{j}\sum_{i=0}^{\gamma-(\inda-\mu-1-j)}\binom{\mu+1}{i}Y^{(\wr_1)}(a,b,\lu;\mu+1+\indb+j-i,i+\inda-\mu-1-j)\lu\nonumber\\
&=\sum_{j=0}^{\nu-\indb}
\binom{\inda-\mu-1}{j}\sum_{i=0}^{\infty}\binom{\inda-\mu-1-j}{i}(-1)^{i}
\big(a_{\inda-j-i}b_{\indb+j+i}\lu
+(-1)^{\inda-\mu-j}b_{\inda+\indb-\mu-1-i}a_{\mu+1+i}
\lu\big)\nonumber\\
&=\sum_{k=0}^{\nu-\indb}\sum_{i=0}^{k}\binom{\inda-\mu-1}{k-i}\binom{\inda-\mu-1-k+i}{i}(-1)^{i}
a_{\inda-k}b_{\indb+k}\lu\nonumber\\
&=\sum_{k=0}^{\nu-\indb}\binom{k-1}{k}a_{\inda-k}b_{\indb+k}\lu\nonumber\\
&=a_{\inda}b_{\indb}\lu.
\end{align}
Writing $\Phi_{\inda,\indb}(a,b,\lu;\wz)=\sum_{i\in\Z}\lambda_{i}\wz^i, \lambda_i\in\C$, by \eqref{eq:sum-Y} we have
\begin{align}
\label{eq:proof-lemma-hom}
\sum_{i\in\Z}\lambda_{i}(a_ib)_{\inda+\indb-i}\lu
&=\sum_{s=0}^{T-1}\sum_{i\in\Z}\lambda_{i}\cty{s}{M}{a,b}{\inda+\indb-i,i}(\lu)\nonumber\\
&=\sum_{s=0}^{T-1}Z_{\inda,\indb}^{(s)}(\Phi_{\inda,\indb}(a,b,\lu;\wz))\nonumber\\
&=Z_{\inda,\indb}^{(\wr_1)}(\Phi_{\inda,\indb}(a,b,\lu;\wz))\nonumber\\
&=a_{\inda}b_{\indb}\lu.
\end{align}
The proof is complete.
\end{proof}

\begin{remark}
\label{remark:(awlb)indau}
The following result, which is complementary to Lemma \ref{lemma:associative-action}, follows from \cite[Lemma 2.8]{Tanabe2015}:
for a weak $(V,T)$-module $\module$, $a,b\in V$, $u\in\module$, and  $\inda\in(1/T)\Z, \wl\in \Z$,
\begin{align}
(a_{\wl}b)_{\inda}u&\in \Span_{\C}\{a_{\wl+\inda-i}b_i\lu\in\module\ |\ i\in(1/T)\Z\}.
\end{align}
Hence a standard argument (cf. \cite[Corollary 4.5.15]{LL}) shows the following result:
let $V$ be a simple vertex algebra, $\module$ a weak $(V,T)$-module, and $a\in V$, $\lu\in\module$. 
If $Y_{\module}(a,x)u=0$, then $a=0$ or $\lu=0$.
\end{remark}

\section{The algebra $\sA_{\alpha}(G,\sS)$ and Proof of Theorem \ref{theorem:main}\label{section:proof}}
Throughout this section, $V$ is a simple vertex algebra of countable dimension and $G$ is a subgroup of  $\Aut V$ of finite order $T$.
In this section, after recalling some properties of the algebra  $\sA_{\alpha}(G,\sS)$ introduced in \cite{DY2002},
we give a proof of Theorem \ref{theorem:main}\label{section:proof}.
\subsection{The algebra $\sA_{\alpha}(G,\sS)$ and its modules}
For a weak $(V,T)$-module $(M,Y_M)$ and $\sigma\in \Aut V$ we define a weak $(V,T)$-module
$(\module\cdot \sigma,Y_{\module\cdot \sigma})$
by $\module\cdot \sigma=\module$ as a vector space and 
\begin{align}
\label{eqn:actionYmodulecdotsigma}
Y_{\module\cdot \sigma}(a,x)&=Y_{\module}(\sigma a,x) \mbox{ for }a\in V.
\end{align}
Note that if $\module$ is irreducible, then so is  $M\cdot \sigma$.
A set $\sS$ of inequivalent irreducible weak $(V,T)$-modules 
is called {\it $G$-stable} if for any $M\in \sS$ and $\sigma \in G$ there exists
$\mN\in \sS$ such that $\module\cdot \sigma \cong \mN$ as weak $(V,T)$-modules.
For a finite right $G$-set $\sX$, a finite dimensional semisimple associative $\C$-algebra $\sA_{\alpha}(G,\sX)$ was constructed  in  \cite{DY2002}.
For  a finite $G$-stable set $\sS$ of inequivalent irreducible weak $(V,T)$-modules,
we recall the definition of  $\sA_{\alpha}(G,\sS)$
and its basic properties according to \cite{DY2002}.
We assume that $\sS$ is a finite $G$-stable set  of inequivalent irreducible weak $(V,T)$-modules until the end of this section. Let
$\module\in\sS$ and $\sigma \in G$. Then there exists $\mN\in\sS$ such that $N\cdot \sigma \cong M$ as 
weak $(V,T)$-modules. That is, there is an isomorphism $\phi(\sigma , \module):\module\rightarrow\mN$ of vector spaces such
that
\begin{align}
\phi(\sigma ,\module)Y_{\module}(a,x)&=Y_{\mN}(\sigma a,x)\phi(\sigma ,\module)
\end{align}
for all $a\in V$. For $\sigma ,\tau\in G$, by the irreducibility of $\module$, there exists $\alpha_{\module}(\sigma ,\tau)\in\C$ such that
\begin{align}
\phi(\sigma ,\module \cdot \tau^{-1})\phi(\tau,\module)=\alpha_{\module}(\sigma ,\tau)\phi(\sigma \tau,\module).
\end{align}
Moreover, for $\sigma ,\tau,\rho\in G$ and $\module\in\sS$ we have
\begin{align}
\alpha_{M\cdot \rho^{-1}}(\sigma,\tau)\alpha_{\module}(\sigma \tau,\rho)&=
\alpha_{\module}(\sigma ,\tau\rho)\alpha_{\module}(\tau,\rho).
\end{align}
We define a vector space $\C\sS=\oplus_{\module\in\sS}\C e(\module)$
with a basis $e(\module)$ for $\module\in\sS$.
The vector space $\C\sS$ is an associative algebra with multiplication given by $e(\module)e(\mN)=\delta_{\module,\mN}e(\mN)$.
We define a subset $\sU(\C \sS)=\{\sum_{\module\in\sS}\lambda_{\module}e(\module)\ |\ \lambda_{\module}\in\C^{\times}\}$ 
of $\C\sS$ where $\C^{\times}=\C\setminus\{0\}$.
The set $\sU(\C\sS)$ becomes a right $G$-set by the action  
\begin{align}
\sum_{\module\in\sS}\lambda_{\module}e(\module)\cdot \sigma &=\sum_{\module\in\sS}\lambda_{\module}e(\module\cdot \sigma )
\end{align}
for $\sigma \in G$.
For $\alpha(\sigma ,\tau)=\sum_{\module\in\sS}\alpha_{\module}(\sigma ,\tau)e(\module)$ we have
\begin{align}
\textcolor{black}{(\alpha(\sigma ,\tau)\cdot \rho)}\alpha(\sigma \tau,\rho)&=\alpha(\sigma ,\tau \rho)\alpha(\tau,\rho)
\end{align}
for all $\sigma ,\tau,\rho\in G$ and hence $\alpha : G\times G\rightarrow \sU(\C \sS)$ is a $2$-cocycle.

The vector space $\sA_{\alpha}(G,\sS)=\C[G]\otimes_{\C}\C \sS$  is an associative $\C$-algebra with multiplication given by
\begin{align}
\label{eq:algebraAalphaGS}
(\sigma\otimes e(\module))\cdot (\tau\otimes e(\mN))&=\alpha_{\mN}(\sigma,\tau)\sigma\tau\otimes e(\module\cdot \tau)e(\mN).
\end{align}
We define an action on $\sA_{\alpha}(G,\sS)$ on $\sM=\oplus_{\module\in \sS}\module$ as follows:
for $\module,\mN\in \sS$, $\lu\in\mN$, and $\sigma\in G$, we define
\begin{align}
\sigma\otimes e(\module)\cdot \lu&=\delta_{\module,\mN}\phi(\sigma,\module)\lu.
\end{align}
Note that the actions of $\sA_{\alpha}(G,\sS)$ and $V^G$ on $\sM$ commute with each other.
For $\module\in \sS$, we define
\begin{align}
G_{\module}&=\{\sigma\in G\ |\ \module\cdot \sigma\cong \module\mbox{ as weak $(V,T)$-modules}\}.
\end{align}  
The restriction map $\alpha_{M}=\alpha|_{G_{M}\times G_{M}}$ is clearly a $2$-cocycle for $G_{M}$.
Let $\sO_{\module}$ be the orbit of $\module$ under the action of $G$
and let $G=\cup_{j=1}^{k}G_{\module}g_{j}$
be a right coset decomposition with $g_1=1$.
Then $\sO_{\module}=\{\module\cdot g_j\ |\ j=1,\ldots,k\}$ and 
$G_{\module\cdot g_j}=g_j^{-1}G_{\module}g_j$.
We define several subspaces of $\sA_{\alpha}(G,\sS)$ by 
\begin{align}
S(\module)&=\Span_{\C}\{\sigma\otimes e(\module)\ |\ \sigma\in G_{\module}\},\nonumber\\
D(\module)&=\Span_{\C}\{\sigma\otimes e(\module)\ |\ \sigma\in G\},\nonumber\\
D(\sO_{\module})&=\Span\{\sigma\otimes e(\module\cdot g_j)\ |\ \sigma\in G, j=1,\ldots,k\}.
\end{align}
The subalgebra $S(\module)$ of $\sA_{\alpha}(G,\sS)$ is isomorphic to 
$\C^{\alpha_{\module}}[G_{\module}]$, the twisted group algebra with $2$-cocycle $\alpha_{\module}$.
Decompose $\sS$ into a disjoint union of orbits:
\begin{align}
\sS = \bigcup_{j\in J}\sO_{j}.
\end{align}
Let $\module^{j}$ be a representative element of $\sO_j$.
Then $\sO_{j}=\sO_{\module^{j}}=\{\module^{j}\cdot\sigma\ |\ \sigma\in G\}$
and 
$\sA_{\alpha}(G,\sS)=\oplus_{j\in J}D(\textcolor{black}{\sO_{\module^{j}}})$.
For $\module\in\sS$, let $\Lambda_{G_{\module}, \alpha_{\module}}$
be the set of all irreducible characters $\lambda$ of $\C^{\alpha_{\module}}[G_{\module}]$.
We denote the corresponding irreducible module by $W(\module)_{\lambda}$. 
Note that $M$ is a semisimple $\C^{\alpha_{\module}}[G_{\module}]$-module. 
Let $M^{\lambda}$ be the sum of irreducible $\C^{\alpha_{\module}}[G_{\module}]$-module
of $\module$ isomorphic to $W(\module)_{\lambda}$. 
Then
\begin{align}
M&=\bigoplus_{\lambda\in\Lambda_{G_{\module}, \alpha_{\module}}}M^{\lambda}.
\end{align}
Moreover $M^{\lambda}=W(\module)_{\lambda}\otimes M_{\lambda}$ where 
\begin{align}
\label{eq:Mlambda=HomC}
M_{\lambda}=\Hom_{\C^{\alpha_{\module}}[G_{\module}]}(W(\module)_{\lambda},\module)
\end{align}
is the multiplicity of $W(\module)_{\lambda}$ in $\module$.
We can realize $M_{\lambda}$ as a subspace of $\module$ in the
following way: let $w\in W(\module)_{\lambda}$ be a fixed nonzero vector. Then we can identify
$\Hom_{\C^{\alpha_{\module}}[G_{\module}]}(W(\module)_{\lambda},\module)$
 with the subspace
\begin{align}
\{f(w)\ |\ f\in \Hom_{\C^{\alpha_{\module}}[G_{\module}]}(W(\module)_{\lambda},\module)\}
\end{align}
of $M^{\lambda}$.
Note that the actions of $\C^{\alpha_{\module}}[G_{\module}]$ and $V^{G_{\module}}$ on $\module$
commute with each other.
So $M^{\lambda}$ and $M_{\lambda}$ are weak $(V^{G_{\module}},T)$-modules.
Furthermore, $\module^{\lambda}$ and $\module_{\lambda}$ are weak 
$(V^G,T)$-modules. The same argument as in \cite[Corollary 6.5 and Theorem 6.7]{DRY2023} shows that
$\module_{\lambda}\neq 0$ and $\module_{\lambda}$ is an irreducible  weak $(V^{G_{\module}},T)$-module
for any $\lambda\in \Lambda_{G_{\module}, \alpha_{\module}}$.

Let $\sS=\cup_{j\in J}\sO_{j}$ be an orbit decomposition and fix $\module^{j}\in \sO_j$ for each $j\in J$.
For convenience, we set
\begin{align}
\label{eq:GjGmoduelj}
G_{j}&=G_{\module^{j}},\quad
\Lambda_{j}=\Lambda_{G_{\module^{j}}, \alpha_{\module^{j}}}, \quad\mbox{ and }\quad
\mW_{j,\lambda}=W(\module^{j})_{\lambda}
\end{align}
for $j\in J$ and $\lambda\in \Lambda_{j}$. We have a decomposition
\begin{align}
M^{j}&=\bigoplus_{\lambda\in\Lambda_{j}}W_{j,\lambda}\otimes M^{j}_{\lambda}
\end{align}
as a $\C^{\alpha_{\module^{j}}}[G_{\module^{j}}]$-module. 
For $j\in J$ and $\lambda\in\Lambda_{j}$ we set
\begin{align}
\label{eq:definitionWjlambda}
\mW_{\lambda}^{j}=
\Ind^{D(\module^{j})}_{S(\module^{j})}W_{j,\lambda}.
\end{align}
It is  shown in \cite[Theorem 3.6]{DY2002} that the algebra $\sA_{\alpha}(G,\sS)$ is semisimple and 
$\{W^{j}_{\lambda}\ |\ j\in J, \lambda\in\Lambda_{j}\}$
is a complete set of representatives of isomorphism
classes of irreducible $\sA_{\alpha}(G,\sS)$-modules.
It follows from \cite[Propositions 6.6 and 6.7]{DY2002} that
under the action of $\sA_{\alpha}(G,\sS)\otimes_{\C} \textcolor{black}{V^{G}}$,
the space $\sM$ decomposes into
\begin{align}
\label{eq:Mdecomposition}
\sM&=\bigoplus_{\lambda\in\Lambda_{j}}\mW_{\lambda}^{j}\otimes M^{j}_{\lambda}.
\end{align}
\subsection{Proof of Theorem \ref{theorem:main}}
We note that \cite[Theorem 4.7]{DRY2023}  holds even for weak $(V,T)$-modules.
Thanks to Lemma \ref{lemma:tensor}, Lemma \ref{lemma:associative-action}, \eqref{eq:Mdecomposition}, 
\cite[Corollary 3.3 and Theorem 4.7]{DRY2023},
the same argument as in \cite[Theorem 7.4]{DRY2023} shows the result.

\providecommand{\MR}{\relax\ifhmode\unskip\space\fi MR }
\providecommand{\MRhref}[2]{%
  \href{http://www.ams.org/mathscinet-getitem?mr=#1}{#2}
}
\providecommand{\href}[2]{#2}

\end{document}